\documentclass[11pt,english]{amsart}
\pdfoutput=1
\usepackage{verbatim}
\usepackage{amsthm}
\usepackage{amstext}
\usepackage{graphicx}

\usepackage[pdfusetitle]{hyperref}
\usepackage{xcolor}
\definecolor{dark-red}{rgb}{0.6,0.15,0.15}
\definecolor{dark-blue}{rgb}{0.15,0.15,0.6}
\definecolor{medium-blue}{rgb}{0,0,0.5}
\hypersetup{
    colorlinks, linkcolor={dark-red},
    citecolor={dark-blue}, urlcolor={medium-blue}
}

\makeatletter

\newcommand{\noun}[1]{\textsc{#1}}
\let\SF@@footnote\footnote
\def\footnote{\ifx\protect\@typeset@protect
    \expandafter\SF@@footnote
  \else
    \expandafter\SF@gobble@opt
  \fi
}
\expandafter\def\csname SF@gobble@opt \endcsname{\@ifnextchar[
  \SF@gobble@twobracket
  \@gobble
}
\edef\SF@gobble@opt{\noexpand\protect
  \expandafter\noexpand\csname SF@gobble@opt \endcsname}
\def\SF@gobble@twobracket[#1]#2{}

\numberwithin{equation}{section}
\numberwithin{figure}{section}
\theoremstyle{plain}
\newtheorem{thm}{\protect\theoremname}
\newenvironment{lyxcode}
{\par\begin{list}{}{
\setlength{\rightmargin}{\leftmargin}
\setlength{\listparindent}{0pt}
\raggedright
\setlength{\itemsep}{0pt}
\setlength{\parsep}{0pt}
\normalfont\ttfamily}%
 \item[]}
{\end{list}}

\renewcommand{\uppercasenonmath}[1]{}

\makeatother

\usepackage{babel}
\providecommand{\theoremname}{Theorem}

\begin{document}

\title{\textup{How Accurate is }\texttt{\textup{inv(A){*}b}}\textup{?}}
\begin{abstract}
Several widely-used textbooks lead the reader to believe that solving
a linear system of equations $Ax=b$ by multiplying $b$ by a computed
inverse \texttt{inv(A)} is inaccurate. Virtually all other textbooks
on numerical analysis and numerical linear algebra advise against
using computed inverses without stating whether this is accurate or
not. In fact, under reasonable assumptions on how the inverse is computed,
\texttt{x=inv(A){*}b} is as accurate as the solution computed by the
best backward-stable solvers. This fact is not new, but obviously
obscure. We review the literature on the accuracy of this  computation
and present a self-contained numerical analysis of it.
\end{abstract}

\author[Alex Druinsky and Sivan Toledo]{Alex Druinsky$^{*}$ and Sivan Toledo$^{*}$}
\thanks{$^{*}$School of Computer Science, Tel Aviv University}

\date{}

\maketitle

\section{Introduction}

Can you accurately compute the solution to a linear equation $Ax=b$
by first computing an approximation $V$ to $A^{-1}$ and then multiplying
$b$ by $V$ (\texttt{x=inv(A){*}b} in Matlab)?

Unfortunately, most of the literature provides a misleading answer
to this question. Many textbooks, including recent and widely-used
ones, mislead the reader to think that \texttt{x=inv(A){*}b} is less
accurate than \texttt{x=A\textbackslash{}b}, which computes the $LU$
factorization of $A$ with partial pivoting and then solves for $x$
using the factors~\emph{\cite[p. 31]{ForsytheMalcolmMoler}}, \emph{\cite[p. 53]{Moler}},
\emph{\cite[p. 50]{Heath}}, \emph{\cite[p. 77]{OLeary}}, \emph{\cite[p. 166]{ConteDeBoor}},
\emph{\cite[pp. 184, 235, and 246]{Stewart}}. Other textbooks warn
against using a computed inverse for performance reasons without saying
anything about accuracy. If you still dare use \texttt{x=inv(A){*}b}
in Matlab code, Matlab's analyzer issues a wrong and misleading warning~\cite{MLint-R2010b}.

As far as we can tell, only two sources in the literature present
a correct analysis of this question. One is almost 50 years old~\cite[pp. 128--129]{Wilkinson},
and is therefore hard to obtain and somewhat hard to read. The other
is recent, but relegates this analysis to the solution of an exercise,
rather than including it in the 27-page chapter on the matrix inverse\emph{~\cite[p. 559; see also p. 260]{Higham}};
even though the analysis there shows that \texttt{x=inv(A){*}b} is
as accurate as \texttt{x=A\textbackslash{}b}, the text ends by stating
that {}``multiplying by an explicit inverse is simply not a good
way to solve a linear system''. The reader must pay careful attention
to the analysis if he or she is to answer our question correctly.

Our aim in this article is to clarify to researchers (and perhaps
also to educators and students) the numerical properties of a solution
to $Ax=b$ that is obtained by multiplying by a computed inverse.
We do not present new results; we present results that are known,
but not as much as they should be. 

Computing the inverse requires more arithmetic operations than computing
an $LU$ factorization. We do not address the question of computational
efficiency, but we do note that there is evidence that using the inverse
is sometimes preferable from the performance perspective~\cite{DitkowskiFibichGavish}. 

It also appears that explicit inverses are sometimes used when the
inverse must be applied in hardware, as in some MIMO radios~\cite{eberli08,StuderEtAl2011}.
The numerical analysis in the literature and in this paper does not
apply as-is to these computations, because hardware implementations
typically use fixed-point arithmetic rather than floating point. Still,
the analysis that we present here provides guiding principles to all
implementations (e.g., to solve for the rows of the inverse using
a backward-stable solver), and it may also provide a template for
an analysis of fixed-point implementations or alternative inversion
algorithms.

The rest of this paper is organized as follows. Section~\ref{sec:A-Loose-Bound}
presents the naive numerical analysis that probably led many authors
to claim that \texttt{x=inv(A){*}b} is inaccurate; the analysis is
correct, but the error bound that it yields is too loose. Section~\ref{sec:Tightening-the-Bound}
presents a much tighter analysis, due to Wilkinson; Higham later showed
that this bound holds even in the componentwise sense. Section~\ref{sec:Left-and-Right}
explains another aspect of computed inverses that is not widely appreciated:
that they are typically good for applying either from the left or
from the right, but not both. Even when \texttt{x=inv(A){*}b} is accurate,
\texttt{x} is usually not backward stable; Section~\ref{sec:backward-stability-of-xv}
discusses conditions under which \texttt{x} is also backward stable.
To help the reader fully understand all of these results, Section~\ref{sec:Numerical-Examples}
demonstrates them using simple numerical experiments. We present concluding
remarks in Section~\ref{sec:Closing-Remarks}.

\section{\label{sec:A-Loose-Bound}A Loose Bound}

Why did the inverse acquire its bad reputation? Good inversion methods
produce a computed inverse $V$ that is, at best, \emph{conditionally}
accurate,
\begin{equation}
\frac{\left\Vert V-A^{-1}\right\Vert }{\left\Vert A^{-1}\right\Vert }=O(\kappa(A)\epsilon_{\text{machine}})\;.\label{eq:conditionally-accurate-V}
\end{equation}
We cannot hope for an unconditional bound of $O(\epsilon_{\text{machine}})$
on the relative forward error. Some inversion methods guarantee conditional
accuracy (for example, computing the inverse column by column using
a backward stable linear solver). In particular, \noun{lapack}'s \texttt{xGETRI}
satisfies (\ref{eq:conditionally-accurate-V}), and also a componentwise
conditional bound~\cite[p. 268]{Higham}. That is, each entry in
the computed inverse that \texttt{xGETRI} produces is conditionally
accurate. It appears that Matlab's \texttt{inv} function also satisfies
(\ref{eq:conditionally-accurate-V}).

Let's try to use (\ref{eq:conditionally-accurate-V}) to obtain a
bound on the forward error $\|x_{V}-x\|$. Multiplying $b$ by $V$
in floating point produces $x_{V}$ that satisfies $x_{V}=\left(V+\Delta\right)b$
for some $\Delta$ with $\|\Delta\|/\|V\|=O(\epsilon_{\text{machine}})$.
Denoting $\Gamma=V-A^{-1}$, we have
\begin{eqnarray*}
x_{V} & = & \left(V+\Delta\right)b\\
 & = & \left(A^{-1}+\Gamma+\Delta\right)b\\
 & = & \left(A^{-1}+\Gamma+\Delta\right)Ax\\
 & = & x+\Gamma Ax+\Delta Ax\;,
\end{eqnarray*}
so
\begin{eqnarray}
\left\Vert x_{V}-x\right\Vert  & \leq & \left\Vert \Gamma\right\Vert \left\Vert A\right\Vert \left\Vert x\right\Vert +\left\Vert \Delta\right\Vert \left\Vert A\right\Vert \left\Vert x\right\Vert \nonumber \\
 & \leq & O(\kappa(A)\epsilon_{\text{machine}})\left\Vert A^{-1}\right\Vert \left\Vert A\right\Vert \left\Vert x\right\Vert +O(\epsilon_{\text{machine}})\left\Vert V\right\Vert \left\Vert A\right\Vert \left\Vert x\right\Vert \nonumber \\
 & = & O(\kappa^{2}(A)\epsilon_{\text{machine}})\left\Vert x\right\Vert +O(\epsilon_{\text{machine}})\left\Vert V\right\Vert \left\Vert A\right\Vert \left\Vert x\right\Vert \;.\label{eq:loose-accuracy-bound-pre}
\end{eqnarray}
Unless $A$ is so ill conditioned that the left-hand side of (\ref{eq:conditionally-accurate-V})
is larger than $1$ (any constant would do), $\|V\|=\Theta(\|A^{-1}\|)$.
Therefore,
\begin{equation}
\left\Vert x_{V}-x\right\Vert \leq O(\kappa^{2}(A)\epsilon_{\text{machine}})\left\Vert x\right\Vert \;.\label{eq:loose-bound}
\end{equation}

In contrast, solving $Ax=b$ using a backward stable solver such as
one based on the $QR$ factorization (or on an $LU$ factorization
with partial pivoting provided there is no growth) yields $x_{\text{backward-stable}}$
for which 
\begin{equation}
\left\Vert x_{\text{backward-stable}}-x\right\Vert \leq O(\kappa(A)\epsilon_{\text{machine}})\left\Vert x\right\Vert \;.\label{eq:backward-stable-accuracy}
\end{equation}

The bound (\ref{eq:loose-bound}) is correct, but it is just an upper
bound on the error, and it turns out that it is loose by a factor
of $\kappa(A)$. It appears that this easy-to-derive but loose bound
gave the matrix inverse its bad reputation. In fact, $x_{V}$ satisfies
an accuracy bound just like (\ref{eq:backward-stable-accuracy}).

\section{\label{sec:Tightening-the-Bound}Tightening the Bound}

The bound (\ref{eq:loose-bound}) is loose because of a single term,
$\left\Vert \Gamma\right\Vert \left\Vert A\right\Vert $, which we
used to bound the norm of $\Gamma A$. The other term in the bound,
$O(\kappa(A)\epsilon_{\text{machine}})\left\Vert x\right\Vert $,
is tight.

The key insight is that rows of $\Gamma=V-A^{-1}$ tend to lie mostly
in the directions of left singular vectors of $A$ that are associated
with small singular values. The smaller the singular value of $A$,
the stronger the influence of the corresponding singular vector (or
singular subspace) on the rows of $\Gamma$. Therefore, the norm of
the product of $\Gamma$ and $A$ is much smaller than the product
of the norms; $A$ shrinks the strong directions of the error matrix
$\Gamma$. This explains why the norm of $\Gamma A$ is small. This
relationship between the singular vectors of $A$ and $\Gamma$ depends
on a backward stability criterion on $V$, which we define and analyze
below.

Suppose that we use a backward stable solver to compute the rows of
$V$ one by one by solving $v_{i}A=e_{i}$ where $e_{i}$ is row $i$
of $I$. Each computed row satisfies 
\[
v_{i}\left(A+\Xi_{i}\right)=e_{i}
\]
with $\|\Xi_{i}\|/\|A\|=O(\epsilon_{\text{machine}})$. Rearranging
the equation, we obtain
\[
VA-I=\left[\begin{array}{c}
-v_{1}\Xi_{1}\\
\vdots\\
-v_{n}\Xi_{n}
\end{array}\right]\;,
\]
so $\|VA-I\|=O(\|V\|\|A\|\epsilon_{\text{machine}})=O(\kappa(A)\epsilon_{\text{machine}})$.
For a componentwise version of this bound and related bounds for other
methods of computing $V$, see~\cite[section 14.3]{Higham}.

This is the key to the conditional accuracy of $x_{V}$. Since $\Gamma A=(V-A^{-1})A=VA-I$,
the norm of $\Gamma A$ is $O(\kappa(A)\epsilon_{\text{machine}})$.
We therefore have the following theorem.
\begin{thm}
\label{thm:inv-accurate}Let $Ax=b$ be a linear system with a coefficient
matrix that satisfies $\kappa(A)\epsilon_{\text{machine}}=O(1)$.
Assume that $V$ is an approximate inverse of $A$ that satisfies
$\|VA-I\|=O(\kappa(A)\epsilon_{\text{machine}})$. Then the floating-point
product $x_{V}$ of $V$ and $b$ satisfies
\[
\frac{\left\Vert x_{V}-x\right\Vert }{\left\Vert x\right\Vert }=O\left(\kappa(A)\epsilon_{\text{machine}}\right)\;.
\]

\end{thm}
The essence of this analysis appears in Wilkinson's 1963 monograph~\cite[pp. 128--129]{Wilkinson}.
Wilkinson did not account for the rounding errors in the multiplication
$Vb$, which are not asymptotically significant, but otherwise his
analysis is complete and correct.

\section{\label{sec:Left-and-Right}Left and Right Inverses}

In this article, we multiply $b$ by the inverse from the left to
solve $Ax=b$. This implies that the approximate inverse $V$ should
be a good left inverse. Indeed, we have seen that a $V$ with a small
left residual $VA-I$ guarantees a conditionally accurate solution
$x_{V}$. Whether $V$ is also a good right inverse, in the sense
that $AV-I$ is small, is irrelevant for solving $Ax=b$. If we were
trying to solve $x^{T}A=b^{T}$, we would need a good right inverse.

Wilkinson noted that if rows of $V$ are computed using $LU$ with
partial pivoting, then $V$ is usually both a good left inverse and
a good right inverse, but not always~\cite[page~113]{Wilkinson}.
Du~Croz and Higham show matrices for which this is not the case,
but they also note that such matrices are the exception rather than
the rule~\cite{DuCrozHigham}.

Other inversion methods tend to produce a matrix that is either a
left inverse or a right inverse but not both. A good example is Newton's
method. If one iterates with $V^{(t)}=(2I-V^{(t-1)}A)V^{(t-1)}$ then
$V^{(t)}$ converges to a left inverse. If one iterates with $V^{(t)}=V^{(t-1)}(2I-AV^{(t-1)})$
then $V^{(t)}$ converges to a right inverse.

Strassen's inversion formula~\cite{BaileyFergusonStrassenInv,Strassen69}
sometimes produces an inverse that is neither a good left inverse
nor a good right inverse~\cite[Section~26.3.2]{Higham}.

\section{\label{sec:backward-stability-of-xv}Multiplication by the Inverse
is (Sometimes) Backward Stable}

The next section presents a simple example in which the computed solution
$x_{V}$ is conditionally accurate but not backward stable. In this
section we show that under certain conditions, the solution is also
backward stable. The analysis also clarifies in what ways backward
stability can be lost.

Suppose that we use a $V$ that is a good right inverse, $\|AV-I\|=O(\kappa(A)\epsilon_{\text{machine}})$.
We can produce such a $V$ by solving for its columns using a backward-stable
solver. We have
\begin{eqnarray*}
Ax_{V}-b & = & A(V+\Delta)b-b\\
 & = & \left(AV-I\right)b+A\Delta b
\end{eqnarray*}
for some $\Delta$ with $\|\Delta\|/\|V\|=O(\epsilon_{\text{machine}})$.
Here too, the $\Delta$ term does not influence the asymptotic upper
bound. The assumption that $V$ is a good right inverse bounds the
other term,
\begin{eqnarray}
\left\Vert Ax_{V}-b\right\Vert  & \leq & \left\Vert AV-I\right\Vert \left\Vert b\right\Vert +\left\Vert A\right\Vert \left\Vert \Delta\right\Vert \left\Vert b\right\Vert \nonumber \\
 & \leq & O(\kappa(A)\epsilon_{\text{machine}})\left\Vert b\right\Vert +O(\epsilon_{\text{machine}})\left\Vert A\right\Vert \left\Vert V\right\Vert \left\Vert b\right\Vert \nonumber \\
 & = & O(\kappa(A)\epsilon_{\text{machine}})\left\Vert b\right\Vert \;.\label{eq:loose-accuracy-bound-pre-1}
\end{eqnarray}
The relative backward error is given by the expression $\|Ax_{V}-b\|/(\|A\|\|x_{V}\|+\|b\|)$~\cite{RigalGaches}.
Filling in the bound on the norm of the residual, we obtain
\begin{eqnarray*}
\frac{\left\Vert Ax_{V}-b\right\Vert }{\|A\|\|x_{V}\|+\|b\|} & \leq & \frac{\left\Vert Ax_{V}-b\right\Vert }{\|A\|\|x_{V}\|}\\
 & = & O\left(\frac{\|A\|\|A^{-1}\|\epsilon_{\text{machine}}\|b\|}{\|A\|\|x_{V}\|}\right)\\
 & = & O\left(\frac{\|A^{-1}\|\|b\|}{\|x_{V}\|}\epsilon_{\text{machine}}\right)\;.
\end{eqnarray*}
If we assume that $V$ is a reasonable enough left inverse so that
at least $\|x_{V}\|$ is close to $\|x\|$ (that is, if the forward
error is $O(1)$), then a solution $x$ that has norm close to $\|A^{-1}\|\|b\|$
guarantees backward stability to within $O(\epsilon_{\text{machine}})$.
Let $A=L\Sigma R^{*}$ be the SVD of $A$, so
\begin{eqnarray*}
x & = & A^{-1}b\\
 & = & R\Sigma^{-1}L^{*}b\\
 & = & \sum_{i}\frac{L_{i}^{*}b}{\sigma_{i}}R_{i}\;,
\end{eqnarray*}
where $L_{i}$ and $R_{i}$ are the left and right singular vectors
of $A$. If $L_{n}^{*}b=O(\|b\|)$, then $\|x\|\geq\left|L_{n}^{*}b\right|/\sigma_{n}=O(\|A^{-1}\|\|b\|)$
and $x_{V}$ is backward stable. If the projection of $b$ on $L_{n}$
is not large but the projection on, say, $L_{n-1}$ is large and $\sigma_{n-1}$
is close to $\sigma_{n}$, the solution is still backward stable,
and so on.

Perhaps more importantly, we have now identified the ways in which
$x_{V}$ can fail to be backward stable:
\begin{enumerate}
\item $V$ is not a good right inverse, or
\item $V$ is such a poor left inverse that $\|x_{V}\|$ is much smaller
than $\|x\|$, or
\item the projection of $b$ on the left singular vectors of $A$ associated
with small singular values is small.
\end{enumerate}
The next section shows an example that satisfies the last condition.

\section{\label{sec:Numerical-Examples}Numerical Examples}

\begin{figure}
\begin{centering}
\includegraphics[width=0.75\textwidth]{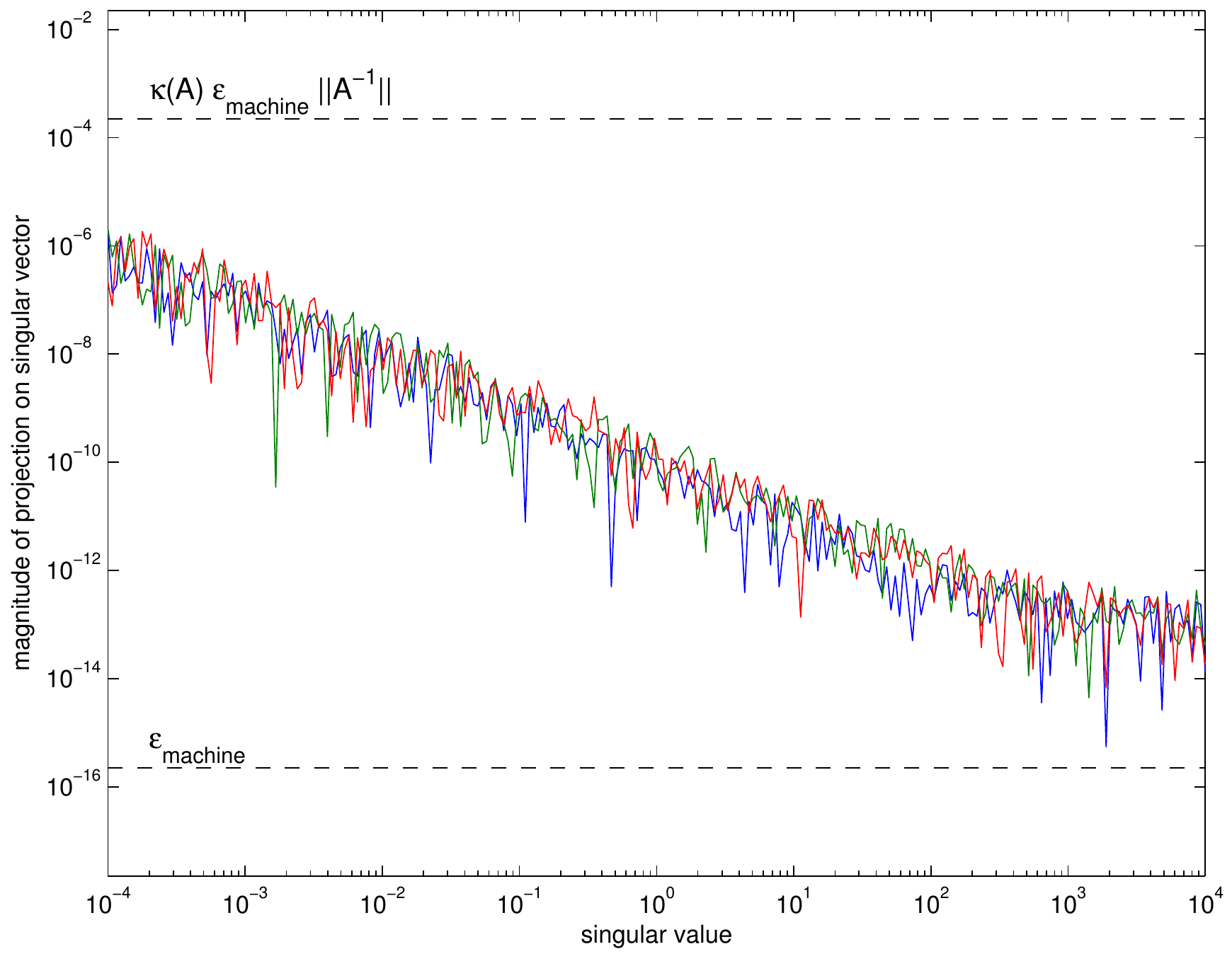}
\par\end{centering}

\caption{\label{fig:proj-Gamma-on-sign-vectors}The magnitude of the projections
of three rows of $\Gamma$ (the first, last, and middle, but all rows
produce similar plots) on the left singular vectors of $A$, as a
function of the corresponding singular values of $A$.}

\end{figure}

Let us demonstrate the theory with a small numerical example. We set
$n=256$, $\sigma_{1}=10^{4}$ and $\sigma_{n}=10^{-4}$, generate
a random matrix $A$ with $\kappa(A)=10^{8}$, and generate its inverse.
The matrix and the inverse are produced by matrix multiplications,
and each multiplication has at least one unitary factor, so both are
accurate to within a relative error of about $\epsilon_{\text{machine}}$.
We also compute an approximate inverse $V$ using \noun{Matlab}'s
\texttt{inv} function.
\begin{lyxcode}
{[}L,dummy,R{]}~=~svd(randn(n));~

svalues~=~logspace(log10(sigma\_1),~log10(sigma\_n),~n);

S~=~diag(svalues);

invS~=~diag(svalues.\textasciicircum{}-1);

A~=~L~{*}~S~{*}~R';

AccurateInv~=~R~{*}~invS~{*}~L';

V~=~inv(A);
\end{lyxcode}
The approximate inverse $V$ is only conditionally accurate, as predicted
by (\ref{eq:conditionally-accurate-V}), but its use as a left inverse
leads to a conditionally small residual.
\begin{lyxcode}
Gamma~=~V~-~AccurateInv;

norm(Gamma)~/~norm(AccurateInv)

~~ans~=~3.4891e-09

norm(V~{*}~A~-~eye(n))

~~ans~=~1.6976e-08
\end{lyxcode}
We now generate a random right hand-side $b$ and use the inverse
to solve $Ax=b$. The result is backward stable to within a relative
error of $\epsilon_{\text{machine}}$. 
\begin{lyxcode}
b~=~randn(n,~1);

x~=~R~{*}~(invS~{*}~(L'~{*}~b));

xv~=~V~{*}~b;

norm(A~{*}~xv~-~b)~/~(norm(A)~{*}~norm(xv)~+~norm(b))

~~ans~=~8.8078e-16~
\end{lyxcode}
Obviously, the solution should be conditionally accurate, and it is.
\begin{lyxcode}
norm(xv~-~x)~/~norm(x)

~~ans~=~3.102e-09~
\end{lyxcode}
We now perform a similar experiment, but with a random $x$, which
leads to a right-hand side $b$ which is nearly orthogonal to the
left singular vectors of $A$ that correspond to small singular values;
now the solution is only conditionally backward stable. 
\begin{lyxcode}
x~=~randn(n,~1);

b~=~L~{*}~(S~{*}~(R'~{*}~x));

xv~=~V~{*}~b;

norm(A~{*}~xv~-~b)~/~(norm(A)~{*}~norm(xv)~+~norm(b))

~~ans~=~2.1352e-10
\end{lyxcode}
Theorem~\ref{thm:inv-accurate} predicts that the solution should
still be conditionally accurate. It is. \noun{Matlab}'s backslash
operator, which is a linear solver based on Gaussian elimination with
partial pivoting, produces a solution with a similar accuracy.
\begin{lyxcode}
norm((A\textbackslash{}b)~-~x)~/~norm(x)

~~ans~=~4.0801e-09

norm(xv~-~x)~/~norm(x)

~~ans~=~4.5699e-09
\end{lyxcode}
The magic is in the special structure of rows of $\Gamma$. Figure~\ref{fig:proj-Gamma-on-sign-vectors}
displays this structure graphically. We can see that a row of $\Gamma$
is almost orthogonal to the left singular vectors of $A$ associated
with large singular values, and that the magnitude of the projections
increases with decreasing singular values. If we produce an approximate
inverse with the same magnitude of error as in \texttt{inv(A)} but
with a random error matrix, it will not solve $Ax=b$ conditionally
accurately.
\begin{lyxcode}
BadInv~=~AccurateInv~+~norm(Gamma)~{*}~randn(n);

xv~=~BadInv~{*}~b;

norm(A~{*}~xv~-~b)~/~(norm(A)~{*}~norm(xv)~+~norm(b))~\\
~~ans~=~0.075727

norm(xv~-~x)~/~norm(x)

~~ans~=~0.83552
\end{lyxcode}

\section{\label{sec:Closing-Remarks}Closing Remarks}

Solving a linear system of equations $Ax=b$ using a computed inverse
$V$ produces a conditionally accurate solution, subject to an easy
to satisfy condition on the computation of $V$. Using Gaussian elimination
with partial pivoting or a $QR$ factorization produces a solution
with errors that have the same order of magnitude as those produced
by $V$.

If the right-hand side $b$ does not have any special relationship
to the left singular subspaces of $A$, then the solution produced
by $V$ is also backward stable (under a slightly different technical
condition on $V$), and hence as good as a solution produced by GEPP
or $QR$. As far as we know, this result is new.

If $b$ is close to orthogonal to the left singular subspaces of $A$
corresponding to small singular values, then the solution produced
by $V$ is conditionally accurate, but usually not backward stable.
Whether this is a significant defect or not depends on the application.
In most applications, it is not a serious problem.

One difficulty with a conditionally-accurate solution that is not
backward stable is that it does not come with a certificate of conditional
accuracy. We normally take a small backward error to be such a certificate.

There might be applications that require a backward stable solution
rather than an accurate one. Strangely, this is exactly the case with
the computation of $V$ itself; the analysis in this paper relies
on rows being computed in a backward-stable way, not on their forward
accuracy. We are not aware of other cases where this is important.

\bibliographystyle{plain}
\bibliography{inv-error}

\begin{thebibliography}{10}

\bibitem{BaileyFergusonStrassenInv}
David~H. Bailey and Helaman R.~P. Ferguson.
\newblock A {S}trassen-{N}ewton algorithm for high-speed parallelizable matrix
  inversion.
\newblock In {\em Proceedings of the 1988 ACM/IEEE conference on
  Supercomputing}, pages 419--424, 1988.

\bibitem{ConteDeBoor}
S.~D. Conte and Carl de~Boor.
\newblock {\em Elementary Numerical Analysis: An Algorithmic Approach}.
\newblock McGraw-Hill Book Company, third edition, 1980.

\bibitem{DitkowskiFibichGavish}
Adi Ditkowski, Gadi Fibich, and Nir Gavish.
\newblock Efficient solution of {$A x^{(k)} = b^{(k)}$} using {$A^{-1}$}.
\newblock {\em Journal of Scientific Computing}, 32:29--44, 2007.

\bibitem{DuCrozHigham}
Jeremy~J. Du~Croz and Nicholas~J. Higham.
\newblock Stability of methods for matrix inversion.
\newblock {\em IMA Journal of Numerical Analysis}, 12(1):1--19, 1992.

\bibitem{eberli08}
Stefan Eberli, Davide Cescato, and Wolfgang Fichtner.
\newblock Divide-and-conquer matrix inversion for linear {MMSE} detection in
  {SDR MIMO} receivers.
\newblock In {\em Proceedings of the 26th Norchip Conference}, pages 162--167.
  IEEE, 2008.

\bibitem{ForsytheMalcolmMoler}
George~E. Forsythe, Michael~A. Malcolm, and Cleve~B. Moler.
\newblock {\em Computer Methods for Mathematical Computations}.
\newblock Prentice-Hall, Inc., 1977.

\bibitem{Heath}
Michael~T. Heath.
\newblock {\em Scientific Computing: An Introductory Survey}.
\newblock The McGraw-Hill Companies, 1997.

\bibitem{Higham}
Nicholas~J. Higham.
\newblock {\em Accuracy and Stability of Numerical Algorithms}.
\newblock SIAM, second edition, 2002.

\bibitem{MLint-R2010b}
The MathWorks, Inc.
\newblock {M}-{L}int: {Matlab}'s code analyzer, 2010.
\newblock Release 2010b.

\bibitem{Moler}
Cleve Moler.
\newblock {\em Numerical Computing with {MATLAB}}.
\newblock SIAM, 2004.

\bibitem{OLeary}
Dianne~P. O'Leary.
\newblock {\em Scientific Computing with Case Studies}.
\newblock SIAM, 2009.

\bibitem{RigalGaches}
J.~L. Rigal and J.~Gaches.
\newblock On the compatibility of a given solution with the data of a linear
  system.
\newblock {\em J. ACM}, 14:543--548, July 1967.

\bibitem{Stewart}
G.~W. Stewart.
\newblock {\em Matrix Algorithms. Volume I: Basic Decompositions}.
\newblock SIAM, 1998.

\bibitem{Strassen69}
V.~Strassen.
\newblock Gaussian elimination is not optimal.
\newblock {\em Numererische Mathematik}, 13:354--355, 1969.

\bibitem{StuderEtAl2011}
C.~Studer, S.~Fateh, and D.~Seethaler.
\newblock {ASIC} implementation of soft-input soft-output {MIMO} detection
  using {MMSE} parallel interference cancellation.
\newblock {\em IEEE Journal of Solid-State Circuits}, 46:1754--1765, 2011.

\bibitem{Wilkinson}
J.~H. Wilkinson.
\newblock {\em Rounding Errors in Algebraic Processes}.
\newblock Prentice Hall, Inc., 1963.

\end{thebibliography}

\end{document}